\def\Theorem#1{{\bf Theorem #1 \qquad }}

\parindent=0 pt  
\parskip = 10 pt  
\baselineskip=15 pt  


Comments on Le Maohua's 1999 paper in the Proc. Japan Acad.  

Reese Scott

revised 26 July 2006

\bigskip

Considering the equation
$$a^x+b^y = c^z \eqno{(1)}$$
where $a$, $b$, $c$ are coprime, squarefree positive integers such that $c$ odd, [Le] gave $2^{\omega(c) + 1}$ as an upper bound on the number of solutions $(x,y,z)$, where $\omega(c)$ is the number of distinct prime factors of $c$.  He also showed that 
$$z < 2 a b \log(2 e a b)/\pi \eqno{(2)}$$
for any solution $(x,y,z)$ to (1).  Here we give slight improvements to each of these results, also removing the restriction that $a$, $b$, $c$ be squarefree:  

\Theorem{}  \quad For positive integers $a$, $b$, $c$, with $c$ odd, (1) has at most $2^{\omega(c)}$ solutions $(x,y,z)$, where $\omega(c)$ is the number of distinct prime factors of $c$.  All solutions $(x,y,z)$ to (1) satisfy $z < a b /2$.  

Proof:  The first assertion follows from the fact that, of the four parity possibilities for the pair $(x,y)$, only two are possible in (1):  this follows from the proof of Theorem 6 of [Sc].  The second assertion follows from Theorem 3 of [Sc-St], noting that $n > n^{1/2} \log(n)$ for $ n \ge 2$.  

\bigskip
  
References 

[Le] M. Le, An upper bound for the number of solutions of the exponential diophantine equation $a^x + b^y = c^z$, Proc. Japan Acad., {\bf 75}, Ser. A (1999)

[Sc]  R. Scott, On the Equations $p^x-b^y = c$ and $a^x+b^y=c^z$, {\it Journal of Number Theory}, {\bf 44}, no. 2 (1993), 153-165.  

[Sc-St]  R Scott and R. Styer, On $p^x - q^y = c$ and related three term exponential Diophantine equations with prime bases, {\it Journal of Number Theory}, {\bf 105} no. 2 (2004), 212--234.

\bye